\documentclass[11pt]{amsart}
\usepackage{amssymb, amscd}
\usepackage{latexsym,epsfig}
\usepackage[all]{xy}
\usepackage{pst-all}
\numberwithin{equation}{section}

\newcommand{\VV}{\mathbb V}

\newcommand{\GG}{\mathbb G}
\newcommand{\PP}{\mathbb P}
\newcommand{\FF}{\mathbb F}

\newcommand{\QQ}{\mathbb Q}
\newcommand{\ZZ}{\mathbb Z}

\newcommand{\A}{\mathcal A}
\newcommand{\M}{\mathcal M}
\newcommand{\s}{\mathbb S}

\newtheorem{theorem}{Theorem}[section]

\newtheorem{corollary}[theorem]{Corollary}
\newtheorem{conjecture}[theorem]{Conjecture}
\newtheorem{definition-lemma}[theorem]{Definition-Lemma}

\theoremstyle{definition}

\newtheorem{notation}[theorem]{Notation}

\theoremstyle{remark}
\newtheorem{remark}[theorem]{Remark}

\begin{document}

\title[Siegel modular forms of level $2$]
{Siegel modular forms of genus 2 and level 2: \\
cohomological computations and conjectures}
\author{Jonas Bergstr\"om}
\address{Korteweg-de Vries Instituut, Universiteit van
Amsterdam, Plantage 
\newline Muidergracht 24, 1018 TV Amsterdam, The Netherlands.}
\email{bergstro@science.uva.nl} 

\author{Carel Faber}
\address{Institutionen f\"or Matematik, Kungliga Tekniska H\"ogskolan,
10044 Stockholm, Sweden.}
\email{faber@math.kth.se}
\author{Gerard van der Geer}
\address{Korteweg-de Vries Instituut, Universiteit van
Amsterdam, Plantage 
\newline Muidergracht 24, 1018 TV Amsterdam, The Netherlands.}
\email{geer@science.uva.nl}

\maketitle
\begin{section}{Introduction}
In this paper we study the cohomology of certain local systems on moduli
spaces of principally polarized abelian surfaces with a level
$2$ structure that corresponds to prescribing a number of Weierstrass
points in case the abelian surface is the Jacobian of a curve of
genus $2$. These moduli spaces are defined over $\ZZ[1/2]$ and 
we can calculate the trace of Frobenius 
on the alternating sum of the \'etale cohomology groups
of these local systems by counting the number of pointed curves of genus $2$ 
with a prescribed number of Weierstrass points that
separately or taken together are defined over the given finite field. 
This cohomology is intimately related to vector-valued Siegel modular forms. 
Two of the present authors carried out this scheme for local 
systems on the moduli space $\A_2$ of level $1$ in \cite{FvdG}.
This provided new information on Siegel modular forms and 
led for example to a precise formulation of a conjecture of Harder 
about congruences between genus $1$ and genus $2$ modular forms and
also to experimental evidence supporting it, cf.\ \cite{Ha,vdG1}.

Here we extend this scheme to level $2$, where new phenomena appear.
In order to be able to extract information on Siegel modular 
forms we must subtract the contributions to the cohomology from the boundary,
that is, the Eisenstein cohomology, and the endoscopic contributions.
We determine the 
contribution of the
Eisenstein cohomology together with its $\s_6$-action
for the full level $2$ structure and on the basis of our computations we 
make precise conjectures on the endoscopic contribution. 
We also make a prediction about the existence of a vector-valued analogue of the
Saito-Kurokawa lift. Assuming these conjectures that are based on ample
numerical evidence, we obtain
the traces of the Hecke-operators $T(p)$ for $p\leq 37$ on the
remaining spaces of `genuine' Siegel modular forms. We present a number
of examples of $1$-dimensional spaces of eigenforms where these traces
coincide with the Hecke eigenvalues to illustrate this.
We hope that the experts on lifting and on endoscopy will be able to prove our
conjectures.
\end{section}
\subsection*{Acknowledgement}
We thank R.\ Tsushima for sending us a program
computing the dimension of the space of Siegel modular cusp forms
of the group $\Gamma_2(w^1)$, used in Section \ref{sec-Tsu}. We also
thank S.\ Kudla for correspondence and
P.\ Deligne for comments on an earlier version.
We also thank S.\ B\"ocherer and R.\ Schulze-Pillot
for pointing out their work on Yoshida-type liftings
and G.\ Harder for discussions on congruences.
C.F.~is supported by the Swedish Research Council (grant 622-2003-1123)
and the G\"oran Gustafsson Foundation for Research in Natural Sciences
and Medicine.
\begin{section}{The Moduli Spaces $\A_2(w^n)$}
Let $\M_2$ be the moduli space of curves of genus $2$
and $\A_2$ the moduli space of principally polarized
abelian surfaces. These are Deligne-Mumford stacks
defined over ${\rm Spec}(\ZZ)$ and by the Torelli map we can view
$\M_2$ as an open substack of $\A_2$.
A curve of genus $2$ admits a unique morphism of degree $2$ to $\PP^1$
and its ramification points are called the Weierstrass points.
If the element $2$ is invertible on the base there are $6$ Weierstrass
points. We now can look at covers of $\M_2$, namely, for
$0\leq n \leq 6$  we consider the stack $\M_2(w^n)$ 
which is the moduli space of pairs $(C,(r_1,\ldots,r_n))$ of curves 
of genus $2$ together with $n$ ordered Weierstrass points. 
These are Deligne-Mumford stacks 
defined over ${\rm Spec}(\ZZ[1/2])$. 

Let $(C,r_1,\ldots,r_6)$ be a genus $2$ curve with its six numbered Weierstrass
points. A Weierstrass point $r_1$ defines an
embedding of $C$ into its Jacobian ${\rm Jac}(C)$ given by $q \mapsto q-r_1$.
This provides us with a set of $5$ points of order $2$ of 
${\rm Jac}(C)$, namely, $\{r_i-r_1: i=2,\ldots,6\}$. The non-zero
points of order $2$ on ${\rm Jac}(C)$ correspond bijectively to
the pairs $\{r_i,r_j\}$ with $i\neq j$.

By associating to a (decorated) curve its
(decorated) Jacobian we have an embedding
$$
\M_2(w^6) \hookrightarrow \A_2[2],
$$ 
where $\A_2[2]$
is the moduli space of principally polarized abelian surfaces $X$ with
a full level $2$-structure, that is, an isomorphism of the kernel
$X[2]$ of multiplication $2_X$ by $2$ on $X$ with a fixed symplectic
module $((\ZZ/2)^4, \langle \, , \, \rangle)$. The symmetric group
$\s_6$ acts on $\M_2(w^6)$ and the group ${\rm GSp}(4,\ZZ/2)$
acts on $\A_2[2]$ and the embedding defines an isomorphism of
$\s_6$ with ${\rm GSp}(4,\ZZ/2)$ which we will make explicit below.
Let us identify ${\rm GSp}(4,\ZZ/2)$ with $\s_6$ under this isomorphism
and define the quotient stacks
$$
\A_2(w^n):= \A_2[2]/\s_{6-n},
$$
where $\s_{6-n}$ is the subgroup of $\s_6$ fixing $\{1,\ldots,n\}$ pointwise.
Note that we have inclusions 
$\M_2(w^n) \hookrightarrow \A_2(w^n)$,
and equalities $\A_2(w^0)=\A_2$ and $\A_2(w^6)=\A_2[2]$.

Let $U$ be a symplectic space of dimension $4$ over $\ZZ/2$
with basis $e_1,e_2,f_1,f_2$ and with symplectic form
$\langle \, , \, \rangle$  such that
$\langle e_i,e_j\rangle=\langle f_i,f_j\rangle =0$ and
$\langle e_i,f_j\rangle= \delta_{i\, j}$ (Kronecker delta).
One observes that $U$ contains six (maximal)
sets $V_i$ ($i=1, \ldots, 6$)
of five vectors $u_j \in U-\{0\}$ with $\langle u_j,u_k\rangle =1$
for $j\neq k$. For instance, one such set is 
$\{e_1,f_1,e_1+f_1+f_2,e_1+e_2+f_1,e_1+e_2+f_1+f_2\}$. 
We have $\# (V_i\cap V_j)=1$ if $i\neq j$.
The action of ${\rm GSp}(4,\ZZ/2)$ on these six sets $V_i$ defines an
isomorphism of ${\rm GSp}(4,\ZZ/2)$ with $\s_6$.
A Weierstrass point $r_i$ on a curve of genus $2$ determines a set
of $5$ points $r_j-r_i$ ($j\neq i$) of order $2$ and they satisfy 
$\langle r_j-r_i,r_k-r_i\rangle=1$ for $j\neq k$. 

Consider the inverse image $\Gamma_2(w^n)$ under
$
{\rm Sp}(4,{\ZZ}) \to
{\rm Sp}(4,{\ZZ}/2) \to \s_6
$
of a subgroup $\s_{6-n}$ of $\s_6$ 
fixing the set $V_i$ for each $i$ between $1$ and $n$. Then the orbifold
$\Gamma_2(w^n)\backslash {\mathcal H}_2$, with ${\mathcal H}_2$
the Siegel upper half space of degree $2$, is the complex fibre
of the moduli stack $\A_2(w^n)$.
By $\Gamma_2[2]=\Gamma_2(w^6)$ we shall mean the kernel of $
{\rm Sp}(4,{\ZZ}) \to {\rm Sp}(4,{\ZZ}/2)$.
\end{section}

\begin{section}{Local Systems}
Let $\pi: {\mathcal X} \to \A_2$ denote the universal
abelian surface. The moduli space $\A_2$ carries natural 
local systems $\VV_{l,m}$
indexed by the pairs $(l,m)$ of integers
with $l\geq m \geq 0$
that correspond to irreducible representations of ${\rm GSp}(4)$,
cf.\ \cite{FvdG}.
The local system $\VV_{1,0}$ is the one defined by $\VV:=R^1\pi_*\QQ_{\ell}$
and $\VV_{l,m}$ is of weight $l+m$ and occurs `for the first time' in
${\rm Sym}^{l-m}(\VV) \otimes {\rm Sym}^m(\wedge^2 \VV)$. When $l > m >0$
the local system is called \emph{regular}.
By pullback to $\A_2(w^n)$
we obtain local systems which we will denote by the same symbol $\VV_{l,m}$.
We are interested in the cohomology (with compact support) of these 
local systems, more precisely in the motivic Euler characteristic
$$
e_c(\A_2(w^n),\VV_{l,m})=
\sum_{i=0}^6 (-1)^i [H_c^i(\A_2(w^n),\VV_{l,m})],
$$
where this expression is taken in the Grothendieck group of an appropriate
category (for instance the category of mixed Hodge structures). Note that
this cohomology is zero if $l+m$ is odd and thus we will from now on only
consider the case when $l+m$ is even.

For any $\A$, among the moduli spaces considered, 
there is a natural map $H^*_c(\A,\VV_{l,m}) \to H^*(\A,\VV_{l,m})$
whose image is called the inner cohomology and denoted by
$H_!^*(\A,\VV_{l,m})$. 
By work of Faltings and Chai one knows that for regular $\lambda$, the 
cohomology groups 
$H_!^i(\A,\VV_{l,m})$ vanish for $i\neq 3$ 
(see \cite{F}, Cor.\ to Thm. 7, p. 84 and
\cite{F-C}, p.~233-7). Moreover, one knows 
that $H^3(\A,\VV_{l,m})$ (resp.\ $H^3_c(\A,\VV_{l,m})$)
carry mixed Hodge structures of weights $\geq l+m+3$ (resp.\ $\leq l+m+3$)
and that $H^3_!(\A,\VV_{l,m})$ carries a pure Hodge structure with Hodge filtration
$$
(0) \subset F^{l+m+3} \subset F^{l+2} \subset F^{m+1} \subset 
F^0=H^3_!(\A,\VV_{l,m}).
$$
The first step in this Hodge filtration is connected to Siegel modular
forms by the isomorphism
$$
F^{l+m+3}\cong S_{l-m,m+3}(\Gamma_2(w^n)),
$$
where $S_{l-m,m+3}(\Gamma_2(w^n))$ is the complex vector 
space of vector-valued cusp forms of
weight $(l-m,m+3)$ on the group $\Gamma_2(w^n)$, cf.\ \cite{F-C}, Thm.\ 5.5,
see also \cite{Getz}, Thm.\ 17. By a Siegel modular form 
of weight $(j,k)$ we mean a vector-valued function on the Siegel upper half
space that transforms with the factor of automorphy:
$$
{\rm Sym}^j(c\tau +d) \det (c\tau +d)^k,
$$
for $(a,b;c,d) \in {\rm Sp}(4,\ZZ)$ or in a subgroup $\Gamma_2(w^n)$.
Thus, the cohomology of the moduli spaces considered is closely connected to 
Siegel modular forms of the corresponding groups. 
In \cite{FvdG} two
of the three present authors studied the cohomology in the case of level $1$
(i.e.\ on $\A_2$) using counts of curves of
genus $2$ of compact type defined over finite fields (cf. also \ \cite{Getz}).

Recall that the Eichler-Shimura theorem says (cf.\ \cite{del})
that for the local system
$\VV_k:={\rm Sym}^k(\VV)$ (with $\VV:=R^1\pi_*\QQ_l$), on the moduli space
$\A_1$ of elliptic curves with universal family
$\pi: {\mathcal X}_1\to \A_1$, one has for even $k\geq 2$,
$$
\begin{aligned}
-e_c(\A_1,\VV_k)=S[{\rm SL}(2,\ZZ),k+2] +1, \\
-e(\A_1,\VV_k)=S[{\rm SL}(2,\ZZ),k+2] + L^{k+1},
\end{aligned}
$$
where we from now on denote by $L:=h^2({\PP}^1)$
the Tate motive of weight $2$ and by 
$S[{\rm SL}(2,\ZZ),k+2]$ the motive of cusp forms of weight $k+2$ 
of ${\rm SL}(2,\ZZ)$ as constructed by Scholl, cf. \cite{Scholl} (see
also \cite{C-F} for an alternative construction). For $k=0$ one can use 
the same formulas if one puts $S[{\rm SL}(2,\ZZ),2]:=-L-1$.

\end{section}
\begin{section}{The Eisenstein cohomology}
The compactly supported cohomology has a natural map to the usual cohomology
and the kernel is called the \emph{Eisenstein cohomology}. 
The corresponding motivic Euler characteristic is denoted by
$
e_{\rm Eis}(\A, \VV_{l,m}).
$
By the {\sl full Eisenstein cohomology} we mean the difference 
between the compactly supported and the usual cohomology, with 
corresponding Euler characteristic,
$$
e_{\rm Eis^f}(\A, \VV_{l,m}):= 
e_c(\A, \VV_{l,m})-e(\A, \VV_{l,m}).
$$
For example, for genus $1$ we have by Eichler-Shimura
$e_{\rm Eis}({\mathcal A}_1,\VV_k)=-1$, and for the full Eisenstein cohomology 
$e_{\rm Eis^f}({\mathcal A}_1,\VV_k)=L^{k+1}-1$.
\begin{remark}\label{Remarkfulleis}
The full Eisenstein cohomology is anti-invariant under Poincar\'e duality
and $e_{\rm Eis}(\A, \VV_{l,m})$ determines the full Eisenstein
cohomology by anti-symmetrizing. The 
converse also holds by considerations of weights  
if $\lambda$ is regular, cf.~\cite{Peters}.
\end{remark}

We shall write $\Gamma(2)$ for the full level $2$ congruence subgroup of ${\rm SL}(2,\ZZ)$ and $\Gamma_0(N)$ for the congruence subgroup of matrices
$(a,b;c,d)$ with $N|c$. For any of these groups $\Gamma$, we will write 
$S_k(\Gamma)$ and $S[\Gamma,k]$
for the space, respectively the motive, of cusp forms of $\Gamma$ 
of weight $k$. The motive $S[\Gamma,k]$ has a Hodge realization that
decomposes as $S_k(\Gamma)\oplus \bar{S}_k(\Gamma)$. 
Let $S_k(\Gamma)^{\rm new}$ denote the subset of newforms
in $S_k(\Gamma)$, and $S[\Gamma,k]^{\rm new}$ the corresponding motive,
constructed in \cite{Scholl}.

\begin{theorem}\label{thm-Eis} 
For regular pairs $(l,m)$ the Eisenstein cohomology of the local system 
$\VV_{l,m}$ on the moduli space $\A_2[2]$ is given by
$$
\begin{aligned}
15 \dim S_{l-m+2}(\Gamma(2)) -15 & \dim S_{l+m+4}(\Gamma(2))\, L^{m+1} +\cr
+ & 15 \begin{cases} S[\Gamma(2),m+2] + 3 & \text{if $m$ even}  \cr
- S[\Gamma(2),l+3]
& \text{if $m$ odd.} \cr
\end{cases}
\cr
\end{aligned}
$$

\end{theorem}
This can be proved as in \cite{vdG2} using the BGG-complex of Faltings-Chai
(see \cite{F-C}), by first computing $e_{\rm Eis^f}(\A_2[2], \VV_{l,m})$ 
and then deducing $e_{\rm Eis}(\A_2[2], \VV_{l,m})$, 
see Remark \ref{Remarkfulleis}.
If the pair $(l,m)$ is not regular we still expect the formula to hold 
as long as we put $S[\Gamma(2),2]:=-L-1$ in case $m=0$  and 
$\dim S_{2}(\Gamma(2)):=-1$ in case $l=m$.

The factors $15$ in the formula come from the fact that the Satake
compactification of $\A_2[2]$ has $15$ one-dimensional and $15$
zero-dimensional boundary components.

The group $\s_6$ acts on $\A_2[2]$ and this induces an action
on the Eisenstein cohomology of a local system. We can decompose this
piece of the cohomology into irreducible representations for $\s_6$.
Note that we can identify $S_k(\Gamma(2))$ with $S_k(\Gamma_0(4))$ 
via the map
$f(z) \mapsto f(2z)$
and the corresponding motive can be split as
$$
S[\Gamma_0(4),k]=S[\Gamma_0(4),k]^{\rm new}+
2\, S[\Gamma_0(2),k]^{\rm new}+3\, S[{\rm SL}(2,\ZZ),k].
$$

We need some notation concerning representations of $\s_6$. Recall
that the irreducible representations are in natural bijective correspondence
with the
partitions of the number $6$. We shall write the representation corresponding 
to the partition $p$ as $s[p]$.  Let us write
$$
\begin{aligned}
A&:=s[3,1^3]+s[2,1^4]={\rm Ind}(s[2]\otimes s[1^4]),\cr
B&:=s[4,2]+s[3,2,1]+s[2^3]={\rm Ind}(s[2]\otimes s[2^2]),\cr
C&:=s[6]+s[5,1]+s[4,2]={\rm Ind}(s[2]\otimes s[4]),\cr
\end{aligned}
$$
where ${\rm Ind}$ denotes induction from $\s_2 \times \s_4 \subset \s_6$
to $\s_6$. We can view $\s_2\times \s_4$ as the stabilizer of a fixed element
in $U-\{0\}$ for the action of $\s_6 \cong {\rm GSp}(4,\ZZ/2)$.
Moreover, the stabilizer of a $2$-dimensional totally isotropic subspace of $U$
corresponds to the stabilizer $H$ of (say) $(12)(34)(56)$ for the
conjugation action of $\s_6$. This group
is isomorphic to the semi-direct product $\s_3 \ltimes (\ZZ/2)^3$. 
The irreducible representations of $H$ are given by pairs of irreducible representations,
one of $(\ZZ/2)^3$, taken up to the action of $\s_3$, and one of its
stabilizer subgroup in $\s_3$.
We put
$$
\begin{aligned}
A'&:=s[4,1^2]+s[3^2]={\rm Ind}_{H}^{\s_6}(s[1^3] \otimes {\rm id}),\cr
B'&:=s[5,1]+s[4,2]+s[3,2,1]={\rm Ind}_{H}^{\s_6}(s[2,1] \otimes {\rm id}), \cr
C'&:=s[6]+s[4,2]+s[2^3]={\rm Ind}_H^{\s_6}(s[3]\otimes {\rm id}).\cr
\end{aligned}
$$
We note that $\dim A = \dim A'=\dim C=\dim C'=15$ while $\dim B=\dim B'=30$.

\begin{notation}
Define $\tau_{N,k}:=\dim S_k(\Gamma_0(N))^{\rm new}$.
\end{notation}

\begin{theorem} \label{thm-EisS6}
For regular pairs $(l,m)$ the contributions in Thm.\ \ref{thm-Eis} to the 
Eisenstein cohomology of the local system $\VV_{l,m}$ can be 
decomposed under the action of $\s_6$ as follows. The term
$15 \dim S_{k}(\Gamma_0(4))$ with $k=l-m+2$ or $k=l+m+4$ decomposes as
$$
\tau_{4,k} \cdot A'  
+ ( \tau_{2,k} + \tau_{1,k}) \cdot B'
+  \tau_{1,k} \cdot C' 
$$
while the term $15 S[\Gamma_0(4),k]$ with $k=m+2$ or $k=l+3$
can be written as
$$
A\otimes S[\Gamma_0(4),k]^{\rm new} +
 B\otimes S[\Gamma_0(2),k]^{\rm new}
+  (B+C)\otimes S[\Gamma_0(1),k] 
$$
and finally $15\cdot 3 \,  L^0= (B+C)\,  L^0$.
\end{theorem}

For $l=m$ we conjecture that $15 \dim S_{2}(\Gamma_0(4))$ decomposes 
as $-C'$, and for $m=0$ that the term $15 S[\Gamma_0(4),2]$ decomposes
as $C \cdot (-L-1)$. The theorem can be proved by the method of
\cite{vdG2} taking into account the action of $\s_6$ on the boundary 
components.

The $\s_6$-decomposition of the Eisenstein cohomology on $\A_2[2]$
allows one to deduce the formulas for $\A_2(w^n)$ for 
$0\leq  n \leq 6$. For example for $\A_2(w^1)$ and 
$\A_2(w^3)$ we find the following. 

\begin{corollary}
For regular $(l,m)$ 
the Eisenstein cohomology of $\VV_{l,m}$ on $\A_2(w^1)$
is given by
$$
\dim S_{l-m+2}(\Gamma_0(2))-\dim S_{l+m+4}(\Gamma_0(2)) \, L^{m+1}
+\begin{cases} 2(S[m+2]+1) & \text{if $m$ even} \cr
-2\, S[l+3] & \text{if $m$ odd.}\cr
\end{cases}
$$
\end{corollary}

\begin{corollary}
For regular $(l,m)$ the Eisenstein cohomology of $\VV_{l,m}$ on 
$\A_2(w^3)$ is given by
$$
\begin{aligned}
4&\dim S_{l-m+2}(\Gamma_0(4)) -4\dim S_{l+m+4}(\Gamma_0(4)) \, L^{m+1}+\cr
+&\begin{cases} 
3S[\Gamma_0(1),m+2]+3S[\Gamma_0(2),m+2]+S[\Gamma_0(4),m+2] +12 &
\text{if $m$ even} \cr
-3S[\Gamma_0(1),l+3]-3S[\Gamma_0(2),l+3]-S[\Gamma_0(4),l+3]  &
\text{if $m$ odd.} \cr
\end{cases}
\end{aligned}
$$
\end{corollary}
\end{section}

\begin{section}{Counting points over finite fields}
Recall that the moduli space $\A_2[2]$ of principally
polarized abelian surfaces with level $2$ structure can be identified
with the moduli space of tuples $(C,r_1,\ldots,r_6)$, 
where $C$ is either an irreducible genus $2$ curve or a 
pair of genus $1$ curves intersecting in one point, and where
$(r_1,\ldots,r_6)$ is a $6$-tuple of marked Weierstrass points.
In the case of two intersecting elliptic curves these are the points of
order $2$ on the two elliptic curves taking the intersection point
as origin on both.

For an odd prime number $p$ we consider this moduli space over the field
$\overline{\FF}_q$, where $q$ is a power of $p$. 
Let $H^i_{\acute{e}t}$ denote the compactly supported 
$\ell$-adic \'etale cohomology. 
The natural action of $\s_6$ on $\A_2[2]$ induces a decomposition of
$H^i_{\acute{e}t}(\A_2[2] \otimes \overline{\FF}_q,\VV_{l,m})$ 
into pieces denoted 
$H^i_{\acute{e}t,\mu}(\A_2[2]\otimes \overline{\FF}_q,\VV_{l,m})$, 
for each irreducible 
representation of $\s_6$ indexed by the partition $\mu$ of $6$. 
We wish to compute the trace of Frobenius $F_q$ on the Euler 
characteristic 
$$e_{\acute{e}t,\mu}(\A_2[2]\otimes \overline{\FF}_q,\VV_{l,m}):=
\sum_i (-1)^i H^i_{\acute{e}t,\mu}(\A_2[2]\otimes \overline{\FF}_q,\VV_{l,m}).
$$
The necessary information to compute this, for any partition $\mu$ and 
pair $(l,m)$, was
found for all odd $q \leq 37$ with the aid of the computer.
We indicate below how this was done.

We will denote by $k$ a finite field 
and by $k_2$ a degree $2$ extension of $k$.

\subsection{Irreducible curves of genus $2$}
Let $P_2(k)\subset k[x]$ be the set of
of all square-free polynomials of degree $5$ or $6$. 
Each element $f \in P_2(k)$ defines a curve $C_f$ of genus $2$
defined by $y^2=f(x)$. 
For each $f \in P_2(k)$ and $k \in \mathcal{K}:=\{\FF_q: 2 \nmid q, q\leq 37\}$
we computed the following: i) the number of points of $C_f$ defined over $k$;
ii) the number of points of $C_f$ defined over $k_2$;  
iii) the fields of definition of all $6$ ramification points of the
canonical map $C_f \to \PP^1$ given by $(x,y) \to x$.

For a partition $\nu$ of $6$ let $P_2(\nu,k)\subset P_2(k)$ be the subset of
polynomials $f$ defining curves $C_f$ which have fields of definition of 
their ramification points given by $\nu$. 
Using the Lefschetz trace formula we can now, for each pair of numbers $n_1,n_2$ and $k \in \mathcal{K}$, compute 
\begin{equation} \label{eq-a}
a(\mathcal{M}_2,\nu,n_1,n_2):=
\sum_{f \in P_2(\nu,k)} a_1(C_f)^{n_1} \cdot a_2(C_f)^{n_2}/|\mathrm{GL}_2(k)|,
\end{equation}
where $a_1(C_f):=\mathrm{Tr}(F_q,H^1_{\acute et}(C_f))$ and $a_2(C_f):=\mathrm{Tr}(F_q^2,H^1_{\acute et}(C_f))$. Note that $|\mathrm{GL}_2(k)|$ 
is the number of $k$-isomorphisms between the curves of $P_2(\nu,k)$. 

\subsection{Pairs of elliptic curves}
Similarly, let $P_1(k)\subset k[x]$ consist of all square-free polynomials
$f(x) \in k[x]$ of degree $3$ and let
$\mathcal{K}':=\{\FF_q,\FF_{q^2}: 2 \nmid q, q\leq 37 \}$
be a collection of finite fields. 
Each element of $P_1(k)$ defines an elliptic curve $C_f$ given by
$y^2=f(x)$ with $x=\infty$ as origin. 
For each element $k \in \mathcal{K}'$ and $f \in P_1(k)$ with corresponding curve $C_f$ we computed the following:
i) the number of points of $C_f$ defined over $k$; ii)
the fields of definition of 
the $3$ affine ramification points
of the map $C_f \to \PP^1$ given by $(x,y) \to x$.

To get the analogue of \eqref{eq-a} for the pairs of elliptic curves 
joined at the origin, we should sum over all possibilities of 
distributing the ramification points and the marked points (which correspond to the monomials $a_1(C_f)^{n_1}a_2(C_f)^{n_2}$) on the two elliptic curves. 
Let us define $a(\A_{1,1},\nu,n_1,n_2)$ to be the sum, over all 
ordered choices of partitions $\rho$ and $\sigma$ of $3$
such that $\nu=\rho+\sigma$ and integers 
$m_1 \leq n_1$ and $m_2 \leq n_2$, of the following. We put $I'_k:=|\mathrm{GL}_2(k)|/(|k|+1)$,
which is the number of $k$-isomorphisms between the curves of $P_1(k)$,
and we divide into two cases
according to if there is an automorphism which interchanges 
the two elliptic curves or not.

{\sl Case i) }
If $\rho \neq \sigma$ or $(m_1,m_2) \neq (n_1-m_1,n_2-m_2)$
we add $1/2$ times
$$
(\sum_{f \in P_1(\rho,k)} a_1(C_f)^{m_1} \cdot a_2(C_f)^{m_2}/I'_k) \cdot (\sum_{f \in P_1(\sigma,k)} a_1(C_f)^{n_1-m_1} \cdot a_2(C_f)^{n_2-m_2}/I'_k).
$$

{\sl Case ii)} 
If $\rho = \sigma$ and $(m_1,m_2) = (n_1-m_1,n_2-m_2)$, 
we have the contribution from pairs of elliptic curves that are defined 
over $k$,
$$
1/2 \cdot \bigl(\sum_{f \in P_1(\rho,k)} a_1(C_f)^{m_1} \cdot a_2(C_f)^{m_2}/I'_k\bigr)^2.
$$
Moreover, if in addition $n_1 = 0$ and $\nu_i = 0$ for all odd $i$, then     
the two elliptic curves together with marked ramification and ordinary points
may also form a conjugate pair. We construct these by taking an 
elliptic curve defined over $k_2$, and  
join it at the origin with its Frobenius conjugate. 
Define the partition $\nu^{1/2}:=[1^{\nu_2} 2^{\nu_4} 3^{\nu_6}]$.
We then add
$$
1/2 \cdot \bigl(\sum_{f \in P_1(\nu^{1/2},k_2)} a_1(C_f)^{n_2}/I'_{k_2}\bigr).
$$
In both these formulas, the factor $1/2$ is due to the extra automorphism. 

\subsection{Adding the contributions from the two strata}
To the irreducible representation of $\s_6$ indexed by the partition $\mu$ of $6$ we can associate $s_{\mu}$, the ordinary Schur polynomial, and to an irreducible representation of the symplectic group $\mathrm{Sp}(4,\QQ)$ indexed by the partition 
$\lambda=[l,m]$, the Schur polynomial $s_{<\lambda>}$, see \cite[Appendix A]{FH}. Written in terms of the power sums $p_i$ we have 
$s_{\mu}=\sum_{\nu} \alpha^{\mu}_{\nu} \cdot p_1^{\nu_1} \cdots p_6^{\nu_6}$ 
and $s_{<\lambda>}=\sum_{n_1,n_2} \beta^{\lambda}_{n_1,n_2} \cdot p_1^{n_1} p_2^{n_2}$ for some rational numbers $\alpha^{\mu}_{\nu}$ and 
$\beta^{\lambda}_{n_1,n_2}$. The trace of Frobenius on $e_{\acute{e}t,\mu}(\A_2[2]\otimes \overline{\FF}_q,\VV_{\lambda})$ is then equal to (compare \cite[Eq. (3.1)]{B2} or \cite[Sec. 4.2]{BT})
\begin{equation*} 
\sum_{n_1,n_2} \sum_{\nu} \alpha^{\mu}_{\nu} \beta^{\lambda}_{n_1,n_2}  \cdot
\bigl(a(\mathcal{M}_2,\nu,n_1,n_2)+a(\A_{1,1},\nu,n_1,n_2)\bigr) \cdot q^{\frac{|\lambda|-n_1-n_2}{2}}.
\end{equation*}

Using these results we have been able to (conjecturally) identify the 
non-Eisenstein pieces of
$e_{c}(\A_2[2],\VV_{l,m})$. 
In this process we have greatly benefited from 
William Stein's tables of modular forms \cite{Stein}. 

\end{section}

\begin{section}{The Strict Endoscopic Part}\label{sec-endoscopy}
The Hecke algebra of ${\rm GSp}(4,\QQ)$ acts on the inner
cohomology of the local system $\VV_{l,m}$, cf.\ e.g. \cite{Ha1},
\cite{F-C}, p.\ 249 and \cite{T}. 
This inner cohomology  on $\A_2(w^n)$ 
also has a Hodge filtration $(0) \subseteq F^{l+m+3} \subseteq F^{l+2} \subseteq 
F^{m+1}\subseteq F^0$. The action of the Hecke algebra respects this
Hodge filtration.
We now consider the irreducible representations $H$ of the Hecke algebra
occurring in the inner cohomology which have the
property that $F^{l+m+3}\cap H=(0)$.
We define the \emph{strict endoscopic part} of the inner cohomology 
of our local system $\VV_{l,m}$ on $\A_2(w^n)$ to be the direct sum
of the representations with that property,
i.e., the part that contains no
contribution from holomorphic (vector-valued) Siegel modular cusp forms
of weight $(j,k)=(l-m,m+3)$. Here we assume that $h^{l+m+3,0}=h^{0,l+m+3}$,
which follows from the fact that eigenforms have totally real eigenvalues
and all representations of ${\rm GSp}(4,\FF_2)$ are defined
over a totally real field, or even $\QQ$.
This endoscopic part should come from the group
${\rm GL}(2,\QQ)\times {\rm GL}(2,\QQ)/\GG_m$.

Assuming that one has a motive $S[\Gamma_2[2],(j,k)]$
associated to the space $S_{j,k}(\Gamma_2[2])$ 
of Siegel modular forms, or more mildly, that 
one can identify the corresponding 
part (of rank $4\dim S_{l-m,m+3}(\Gamma_2[2])$) in the inner cohomology 
of $\VV_{l,m}$ on $\A_2[2]$, 
one might define the strict endoscopic part of the inner cohomology by 
the equation
$$
e_{\rm End^s}(\A, \VV_{l,m}):= 
e_c(\A, \VV_{l,m})-e_{\rm Eis}(\A, \VV_{l,m})+S[\Gamma_2[2],(l-m,m+3)].
$$

\begin{remark} In the case of a Saito-Kurokawa lift of weight $(0,k)$
the associated rank $4$ motive is $S[SL(2,\ZZ),2k-2]+L^{k-1}+L^{k-2}$.
\end{remark}

Using our counts of curves over finite fields with marked Weierstrass points and
our formulas for the Eisenstein cohomology, we are able to
make precise conjectures on the strict endoscopic part. This part of the cohomology
is describable in terms of elliptic modular forms (but note that there are 
also other parts that are describable in terms of 
such forms).
\begin{conjecture} \label{conj-strEnd}
Let $k:=l+m+4$, $k':=l-m+2$, then the strict endoscopic part of the inner cohomology of $\VV_{l,m}$ on $\A_2[2]$ is given by
$$
 -  5 L^{m+1} \cdot \dim S_k(\Gamma_0(4)) \cdot S[\Gamma_0(4),k'],
$$
where we interpret $S[\Gamma_0(4),2]$ as $-L-1$.
\end{conjecture}

\end{section} 

\begin{section}{A lifting to Vector-valued Modular Forms}\label{sec-lift}
The Saito-Kurokawa lifting (see e.g.\ \cite{K, Z})
gives a way to transform a cusp form $f$ 
that is a normalized eigenform of weight $2k$ ($k$ odd) on 
${\rm SL}(2,\ZZ)$ into 
a scalar-valued cusp (eigen) form of weight $k+1$ on ${\rm Sp}(4,\ZZ)$. 
In terms of Galois representations or
$L$-factors, the reciprocal of the characteristic polynomial of Frobenius 
at a prime $p$ is in the Saito-Kurokawa case equal to 
$$
(1-p^{k-1}X)(1-a(p)X+p^{2k-1}X^2)(1-p^{k}X),
$$
with $a(p)$ the Hecke eigenvalue of $f$ at $p$. 

Based on our calculations of the cohomology of local systems $\VV_{l,m}$
on $\A_2[2]$, we conjecture the following (Yoshida type) 
lifting from pairs
of elliptic modular forms to vector-valued Siegel modular forms.
Recall the notion of spinor $L$-function (see \cite{An}) 
and that the Atkin-Lehner involution $w_2$
acts on $S_k(\Gamma_0(2))$ with eigenspaces $S_k^+(\Gamma_0(2))$ and
$S_k^{-}(\Gamma_0(2))$ for the eigenvalues $+1$ and $-1$.

\begin{conjecture}\label{conj-lift} 
For an eigenform $f \in S_{l+m+4}(\Gamma_0(2))^{\rm new}$ and 
an eigenform $g\in S_{l-m+2}(\Gamma_0(2))^{\rm new}$ there is a Siegel
modular form $F\in S_{l-m,m+3}(\Gamma_2[2])$, an eigenform for the Hecke 
algebra, with spinor $L$-function
$$
L(F,s)=L(f,s)L(g,s-m-1).
$$
The form $F$ will appear with multiplicity $5$ in $S_{l-m,m+3}(\Gamma_2[2])$
if $f$ and $g$ have the same eigenvalue $\pm$ under $w_2$
and with multiplicity $1$ if they have opposite eigenvalues under $w_2$.

Similarly, for an eigenform $f\in S_{l+m+4}(\Gamma_0(4))^{\rm new}$ and
an eigenform $g\in S_{l-m+2}(\Gamma_0(4))^{\rm new}$ there is a Siegel
modular form $F\in S_{l-m,m+3}(\Gamma_2[2])$ with spinor $L$-function
$$
L(F,s)=L(f,s)L(g,s-m-1)
$$ 
and  it will appear with multiplicity $5$ in 
$S_{l-m,m+3}(\Gamma_2[2])$. 
\end{conjecture}

\begin{remark} The first part of this conjecture is consistent with work of
B\"ocherer and Schulze-Pillot who constructed a Yoshida-type
lifting for a pair $(f,g)$ of newforms on $\Gamma_0(2)$ with the
same sign under $w_2$ to Siegel modular forms on the congruence
subgroup $\Gamma_0^{(2)}(2) \subset {\rm Sp}(4,\ZZ)$; see \cite{B-S}, Thm.\ 5.1
and the ensuing remark on p.\ 99.
\end{remark}

Note that by Tsushima (see \cite{Tsushima}) we know the dimensions
of the spaces $S_{j,k}(\Gamma_2[2])$ of Siegel modular forms of weight $(j,k)$
on the group $\Gamma_2[2]$. In the cases $(j,k)=(4,4)$, $(6,3)$ and $(8,3)$
it seems that the
conjectured lifts generate all of $S_{j,k}(\Gamma_2[2])$.

Let us define $S^{\rm lift}_{j,k}(\Gamma_2[2])$ to be the subspace in 
$S_{j,k}(\Gamma_2[2])$ consisting of the cusp forms obtained by the 
lifting described above. The following conjecture tells us 
the action of $\s_6$ on the space of lifted cusp forms.

\begin{notation}
Put
$\tau^{+}_{k}:=\dim S^{+}_k(\Gamma_0(2))^{\rm new}$ and 
$\tau^{-}_{k}:=\dim S^{-}_k(\Gamma_0(2))^{\rm new}$.
\end{notation}

\begin{conjecture} \label{conj-lift-symm}
If we assume that $l \neq m$ and let $k:=l+m+4$, $k':=l-m+2$,  
then $S^{\rm lift}_{l-m,m+3}(\Gamma_2[2])$ decomposes as a representation of
$\s_6$ as
\begin{multline*}
 s[2,1^4]^{\oplus \tau_{4,k'}}  \otimes S_k(\Gamma_0(4))^{\rm new} \oplus 
(s[2^3]^{\oplus \tau^+_{k'}} \oplus s[1^6]^{\oplus \tau^-_{k'}}) \otimes S^{+}_k(\Gamma_0(2))^{\rm new} \oplus \\
(s[2^3]^{\oplus \tau^-_{k'}} \oplus s[1^6]^{\oplus \tau^+_{k'} }) \otimes S^{-}_k(\Gamma_0(2))^{\rm new}.
\end{multline*}
\end{conjecture}

\begin{remark} Note that there are no vector-valued lifts of level $1$. 
\end{remark}

We also give a corresponding conjecture for the Saito-Kurokawa lifts.

\begin{conjecture} \label{conj-SK-symm}
If we assume that $l = m$ and let $k:=l+m+4$, then 
$S^{\rm lift}_{l-m,m+3}(\Gamma_2[2])$ decomposes as a representation of
$\s_6$ for $m$ odd as
$$
s[4,2] \otimes S^+_k(\Gamma_0(2))^{\rm new} \oplus s[2^3] \otimes S^-_k(\Gamma_0(2))^{\rm new}
\oplus (s[6] \oplus s[4,2] \oplus s[2^3]) \otimes S_k(\Gamma_0(1)),
$$
and for $m$ even as
$$
s[3^2] \otimes S_k(\Gamma_0(4))^{\rm new} \oplus s[5,1] \otimes S^-_k(\Gamma_0(2))^{\rm new} \oplus s[1^6] \otimes S^+_k(\Gamma_0(2))^{\rm new}.
$$
\end{conjecture}

\end{section}
\begin{section}{Decomposing the endoscopic contribution}
The Saito-Kurokawa lift for level $1$ 
associates, for odd $l=m$, to the space $S_{l+m+4}({\rm SL}(2,\ZZ))$ 
of cusp forms 
on ${\rm SL}(2,\ZZ)$ the motive 
$$-S[{\rm SL}(2,\ZZ),l+m+4]-s_{l+m+4}(L^{l+2}+L^{m+1})
$$ 
in the cohomology of the local system $\VV_{l,m}$ on ${\mathcal A}_2$, 
the minus sign indicating that it lands in odd degree cohomology. 
In Conjecture 4.1 of
\cite{FvdG} we conjecture the existence of a (strict) endoscopic part
$$
-s_{l+m+4}S[{\rm SL}(2,\ZZ),l-m+2]L^{m+1}=s_{l+m+4}(L^{l+2}+L^{m+1}).
$$ 
Assuming this and adding the two contributions the net result would be
the existence of $S[{\rm SL}(2,\ZZ),l+m+4]$ in the inner cohomology.
In level $2$ we see a similar phenomenon that becomes clearer if one
takes the action of $\s_6$ into account. So it makes sense to add the two
contributions under the heading `expanded endoscopy.'

For $l \neq m$ and $f \in  S_{l+m+4}(\Gamma_0(N))^{\rm new}$ 
and $g \in S_{l-m+2}(\Gamma_0(N))^{\rm new}$ the motive of the 
corresponding lifting is of the form 
$M_f+L^{m+1}M_g$ where $M_f$ and $M_g$
denote the motives associated to the cusp forms $f$ and $g$.
Let us call $M_f$ the `leading' part of the vector-valued lift
and $L^{m+1}M_g$ the `trailing' one. Note that in the cohomology where this lift
appears, the trailing part is of precisely 
the same form as the strict endoscopy. Let us define the 
\emph{expanded endoscopic part} of the cohomology to be the strict 
endoscopy plus the contribution from the trailing terms coming from
the lifts described in Conjecture \ref{conj-lift}.

\begin{conjecture}\label{conj-endoscopy}
Assume that $l \neq m$ and let $k:=l+m+4$, $k':=l-m+2$. The expanded 
endoscopic part
 of the inner cohomology of $\VV_{l,m}$ on $\A_2[2]$ is given by
\begin{multline*}
 - L^{m+1} \Bigl( \bigl(\tau_{4,k} \cdot s[3,1^3]+\tau_{1,k} \cdot s[3^2]+(\tau_
{1,k}+\tau_{2,k}) \cdot s[4,1^2] \bigr) S[\Gamma_0(4),k']^{\rm new} \\
+\bigl((\tau_{1,k}+\tau_{2,k}) \cdot s[3,2,1]+\tau_{4,k} \cdot s[4,1^2] + \tau_{
1,k} \cdot s[4,2]+\tau_{1,k} \cdot s[5,1] \bigr) S[\Gamma_0(2),k']^{\rm new} \\
+\bigl(\tau^+_{k} \cdot s[4,2]+\tau^-_{k}\cdot s[5,1]\bigr) S^+[\Gamma_0(2),k']^
{\rm new}\\
+\bigl(\tau^-_{k} \cdot s[4,2]+\tau^+_{k} \cdot s[5,1]\bigr) S^-[\Gamma_0(2),k']
^{\rm new}\\
 +\bigl(\tau_{1,k} \cdot s[2^3]+(\tau_{1,k}+\tau_{2,k}) \cdot s[3,2,1]+\tau_{4,k} \cdot s[3^2]+\tau_{4,k} \cdot s[4,1^2]+(\tau_{2,k}+2\tau_{1,k})\cdot s[4,2]\\+ (\tau_{1,k}+\tau_{2,k}) \cdot s[5,1]+\tau_{1,k} \cdot s[6]\bigr) S[\Gamma_0(1),k']^{\rm new}
\Bigr).
\end{multline*}
\end{conjecture}

\begin{remark} In Conjecture \ref{conj-endoscopy} 
both the strict endoscopy and the lifts from Conjecture \ref{conj-lift}
contribute to the terms
$-L^{m+1} \cdot \tau_{4,k} \cdot s[3,1^3] S[\Gamma_0(4),k']^{\rm new}$, 
$-L^{m+1} \bigl(\tau^+_{k} \cdot s[4,2]+\tau^-_{k}\cdot s[5,1]\bigr)
S^+[\Gamma_0(2),k']^{\rm new}$ and
$-L^{m+1} \bigl(\tau^-_{k} \cdot s[4,2]+\tau^+_{k} \cdot s[5,1]\bigr)
S^-[\Gamma_0(2),k']^{\rm new}$.
\end{remark}

\begin{conjecture}\label{conj-nonreg-endoscopy}
Assume that $l=m$ and let $k:=2m+4$. The expanded endoscopic part
 of the inner cohomology of $\VV_{m,m}$ on $\A_2[2]$ is given by
\begin{equation*}
L^{m+1} \cdot (L+1) \cdot 
\begin{cases}
\tau^{+}_{k} \cdot s[1^6]+\tau_{4,k} \cdot s[3^2]+\tau^{-}_{k} \cdot s[5,1] & \text{if $m$ odd}  \\
(\tau^{-}_{k}+\tau_{1,k}) \cdot s[2^3]+(\tau^{+}_{k}+\tau_{1,k})\cdot s[4,2]+\tau_{1,k} \cdot s[6] & \text{if $m$ even}
\end{cases}
\end{equation*}
\end{conjecture}

\end{section}

\begin{section}{Dimension Checks} \label{sec-Tsu}
In the case of one Weierstrass point we have computed the numerical
Euler characteristic
$\sum (-1)^i \dim H^i_c(\A_2(w^1),\VV_{l,m}) \in \ZZ$ for any $(l,m)$
using methods as in \cite{JBvdG} and \cite{BvdG}. The conjectural results
agree for $(l,m)$ with $l+m\leq 10$ with these numerical
Euler characteristics of the local systems  and moreover for larger
values of $(l,m)$, e.g.\ for $l+m\leq 20$, we find that the numerical Euler
characteristic minus the conjectured Eisenstein and
endoscopic part, is always a non-positive multiple of $4$.
For all $l+m\leq 20$, this number equals $-4$ times the dimension
of the space of Siegel
modular cusp forms $S_{l-m,m+3}(\Gamma_2(w^1))$ as calculated
by a program provided to us
by R.\ Tsushima.
\end{section}

\begin{section}{Examples of Eigenvalues of Hecke Eigenforms}
We will now give a number of examples.
We first write out some (conjectural) results for the first few local
systems. Needless to say they are based on ample
numerical evidence. Recall that the cohomology has the following parts,
$$e_c(\A, \VV_{l,m})=e_{\rm Eis}(\A, \VV_{l,m})+e_{\rm End^s}(\A, \VV_{l,m})
-S[\Gamma_2[2],(l-m,m+3)].
$$
Here we will write $\Phi_{N,k}:=S[\Gamma_0(N),k]^{\rm new}$ and in all of the
following cases this will be a motive associated to a single newform.

\bigskip
\vbox{
\bigskip\centerline{\def\quad{\hskip 0.6em\relax}
\def\quod{\hskip 0.5em\relax }
\vbox{\offinterlineskip
\hrule
\halign{&\vrule#&\strut\quod\hfil#\quad\cr
height2pt&\omit&&\omit&\cr
& $(l,m)$ && $e_c({\mathcal A}_2[2],\VV_{l,m})$  &\cr
\noalign{\hrule}
height2pt&\omit&&\omit&\cr
  & $(0,0)$ &&  $L^3+L^2-14L+16$ & \cr
  & $(2,0)$ &&  $-30L+30$ & \cr
  & $(1,1)$ &&  $5L^3-10L^2$ & \cr
  & $(4,0)$ &&  $-45 L + 45-10 L \Phi_{4,6}$ & \cr
  & $(3,1)$ &&  $-30 L^2 -15 \Phi_{4,6}$ & \cr
  & $(2,2)$ &&  $9L^4 -21 L^3- \Phi_{2,8}$ & \cr
  & $(6,0)$ &&  $-60L+60 -31L \Phi_{2,8}- \Phi_{2,10}$ & \cr
  & $(5,1)$ &&  $-45 L^2+15 -30\Phi_{2,8}- 20L \Phi_{4,6}-5 \Phi_{4,10}$ 
& \cr
  & $(4,2)$ &&  $-45 L^3+45-S[\Gamma_2[2],(2,5)]$ & \cr
  & $(3,3)$ &&  $10 L^5-35L^4-15\Phi_{4,6}- 5\Phi_{2,10}$ & \cr
  & $(8,0)$ &&  $-75L+75-25L\Phi_{4,10}-40L\Phi_{2,10}-5\Phi_{4,12}$ & \cr
  & $(7,1)$ &&  
$-60L^2+30-15\Phi_{4,10}-30\Phi_{2,10}-40L^2\Phi_{2,8}-S[\Gamma_2[2],(6,4)]$ 
& \cr
  & $(6,2)$ &&  $-60L^3+60-20L^3\Phi_{4,6}-S[\Gamma_2[2],(4,5)]$ & \cr
  & $(5,3)$ &&  $-60 L^4 -30\Phi_{2,8}-S[\Gamma_2[2],(2,6)]$ & \cr
  & $(4,4)$ &&  $15L^6-45L^5+30-15\Phi_{4,6}-5\Phi_{4,12}$ & \cr
} \hrule}
}}
\bigskip

In a number of cases, we can identify ``genuine'' Siegel modular forms,
i.e.\ not lifts of the type described in Conjecture \ref{conj-lift}. 
The space $S_{j,k}(\Gamma_2[2])$, and therefore
also the motive (or corresponding part of the cohomology) 
$S[\Gamma_2[2],(j,k)]$ can be
decomposed under the action of $\s_6$ into a sum of
spaces $S_{j,k}(\Gamma_2[2])^{\mu}$ corresponding to
the $\s_6$-representation given by the partition $\mu$.
For the cases appearing in the table above we have
\begin{gather*}
S_{2,5}(\Gamma_2[2])=S_{2,5}(\Gamma_2[2])^{[2^2,1^2]}\\
S_{6,4}(\Gamma_2[2])=S_{6,4}(\Gamma_2[2])^{[2^2,1^2]} \oplus 
S_{6,4}(\Gamma_2[2])^{[3,1^3]}\\
S_{4,5}(\Gamma_2[2])=S_{4,5}(\Gamma_2[2])^{[2,1^4]} \oplus 
S_{4,5}(\Gamma_2[2])^{[2^2,1^2]} \oplus S_{4,5}(\Gamma_2[2])^{[3,2,1]} \\
S_{2,6}(\Gamma_2[2])=S_{2,6}(\Gamma_2[2])^{[3,1^3]} \oplus 
S_{2,6}(\Gamma_2[2])^{[3,2,1]},
\end{gather*}
and each of these subspaces is generated by one vector-valued Siegel 
modular form. For instance for $(l,m)=(4,2)$ we have one vector-valued 
Siegel modular form appearing with the representation $s[2^2,1^2]$, i.e. 
with multiplicity $9$, which agrees with the result of Tsushima that
$S_{2,5}(\Gamma_2[2])$ is $9$-dimensional, see \cite{Tsushima}. Moreover,
according to our data
$S_{4,5}(\Gamma_2[2])^{[2,1^4]}$ is generated by a lift with corresponding
motive $\Phi_{4,12}+L^3\Phi_{4,6}$.

The trace of Frobenius, for a prime $p>2$, on the space
$S[\Gamma_2[2],(j,k)]$ is
equal to the trace of the Hecke operator $T(p)$ on $S_{j,k}(\Gamma_2[2])$.
In the following two tables we write the (conjectural) Hecke eigenvalues
for the generating Siegel modular form  for $3\leq p \leq 19$
in four cases when $S_{j,k}(\Gamma_2[2])^{\mu}$ is generated by a
single vector-valued Siegel modular form. We are assuming here
the conjectures on the endoscopy given above.
Note that all these eigenvalues have many small prime factors.

\bigskip
\centerline{
\vbox{
\offinterlineskip
\hrule
\halign{&\vrule#& \quad \hfil#\hfil \strut \quad \cr
height2pt&\omit&&\omit&&\omit& \cr
& $p$  && $S_{2,5}(\Gamma_2[2])^{[2^2,1^2]}$ && 
$S_{6,4}(\Gamma_2[2])^{[2^2,1^2]}$  & \cr
height2pt&\omit&&\omit&&\omit&\cr
\noalign{\hrule}
height2pt&\omit&&\omit&&\omit&\cr
& $3$  && $-2^3 \cdot  5$ && $-2^3\cdot 5\cdot 7$    &\cr
height2pt&\omit&&\omit&&\omit&\cr
& $5$  && $-2^2\cdot 5^2\cdot 13$ && $-2^2\cdot 5\cdot 149$   &\cr
height2pt&\omit&&\omit&&\omit&\cr
& $7$  && $2^4\cdot 3\cdot 5\cdot 13$ &&  $-2^4\cdot 3\cdot 5\cdot 401$  
&\cr
height2pt&\omit&&\omit&&\omit&\cr
& $11$  && $2^3\cdot 11\cdot 13\cdot 31$ && $2^3\cdot 36383$   &\cr
height2pt&\omit&&\omit&&\omit&\cr
& $13$  &&  $-2^2\cdot 5\cdot 3469$  && $2^2\cdot 5\cdot 37\cdot 251$  &\cr
height2pt&\omit&&\omit&&\omit&\cr
& $15$  && $-2^2\cdot 5\cdot 11\cdot 13\cdot 197$ &&  $2^2\cdot 5\cdot 
19\cdot 6983$  &\cr
height2pt&\omit&&\omit&&\omit& \cr
& $17$  && $2^3\cdot 5^2\cdot 11\cdot 13\cdot 17$ && $-2^3\cdot 5\cdot 
29\cdot 6287$  &\cr
height2pt&\omit&&\omit&&\omit&\cr
& $19$  && $-2^4\cdot 5\cdot 13\cdot 311$  &&  $-2^4\cdot 5\cdot 43\cdot 
2267$  &\cr
height2pt&\omit&&\omit&&\omit&\cr
} \hrule}}

\bigskip
\centerline{
\vbox{
\offinterlineskip
\hrule
\halign{&\vrule#& \quad \hfil#\hfil \strut \quad \cr
height2pt&\omit&&\omit&&\omit& \cr
& $p$  && $S_{6,4}(\Gamma_2[2])^{[3,1^3]}$ && 
$S_{10,3}(\Gamma_2[2])^{[2^2,1^2]}$  & \cr
height2pt&\omit&&\omit&&\omit&\cr
\noalign{\hrule}
height2pt&\omit&&\omit&&\omit&\cr
& $3$  && $-2^3\cdot 3$  && $2^3\cdot 5^2$    &\cr
height2pt&\omit&&\omit&&\omit&\cr
& $5$  && $2^2\cdot 3^2\cdot 7\cdot 41$ && $2^2\cdot 5\cdot 127$   &\cr
height2pt&\omit&&\omit&&\omit&\cr
& $7$  && $2^4\cdot 5^2\cdot 73$ && $-2^4\cdot 3\cdot 5^2\cdot 13$  &\cr
height2pt&\omit&&\omit&&\omit&\cr
& $11$  && $-2^3\cdot 3^2\cdot 4793$ && $-2^3\cdot 439\cdot 1123$  &\cr
height2pt&\omit&&\omit&&\omit&\cr
& $13$  && $-2^2\cdot 7\cdot 21563$  && $2^2\cdot 5^2\cdot 47\cdot 
4457$  &\cr
height2pt&\omit&&\omit&&\omit&\cr
& $15$  && $-2^2\cdot 3^2\cdot 2351$ && $2^2\cdot 5^2\cdot 799441$  &\cr
height2pt&\omit&&\omit&&\omit& \cr
& $17$  && $-2^3\cdot 7\cdot 11\cdot 37\cdot 383$ && $2^3\cdot 5\cdot 
7\cdot 461\cdot 1723$ &\cr
height2pt&\omit&&\omit&&\omit&\cr
& $19$  && $-2^4\cdot 3^2\cdot 11\cdot 17\cdot 29\cdot 43$ && $2^4\cdot 
5^2\cdot 3653483$  &\cr
height2pt&\omit&&\omit&&\omit&\cr
} \hrule}}
\bigskip

We compute the slopes for the single Siegel modular cusp form generating 
the space $S_{2,6}(\Gamma_2[2])^{[3,1^3]}$ and for the one
generating  $S_{2,6}(\Gamma_2[2])^{[3,2,1]}$.
Recall that the reciprocal of the characteristic polynomial of Frobenius is
$$
1-\lambda(p)X+(\lambda(p)^2-\lambda(p^2)-p^{l+m+2})X^2-\lambda(p)p^{l+m+3}X^3
+p^{2l+2m+6}X^4
$$
and the `slope' refers to the slopes of the Newton polygon.

\bigskip
\centerline{
\vbox{
\offinterlineskip
\hrule
\halign{&\vrule#& \quad \hfil#\hfil \strut \quad \cr
height2pt&\omit&&\omit&&\omit&&\omit& \cr
& $p$  && $\lambda(p)$ && $\lambda(p^2)$  && slopes & \cr
height2pt&\omit&&\omit&&\omit&&\omit&\cr
\noalign{\hrule}
height2pt&\omit&&\omit&&\omit&&\omit&\cr
& $3$  && $2^3\cdot 3^3$  && $-2^2\cdot 3^6\cdot 107$ && $3,3,8,8$ &\cr
height2pt&\omit&&\omit&&\omit&&\omit&\cr
& $5$  && $-2^2\cdot 3^4\cdot 17$ && $2^2\cdot 181\cdot 26161$ && 
$11/2$  &\cr
height2pt&\omit&&\omit&&\omit&&\omit&\cr
\noalign{\hrule}
height2pt&\omit&&\omit&&\omit&&\omit&\cr
\noalign{\hrule}
height2pt&\omit&&\omit&&\omit&&\omit&\cr
height2pt&\omit&&\omit&&\omit&&\omit&\cr
& $3$  && $-2^3\cdot 3^2\cdot 5$  && $3^4\cdot 1753$ && $2,2,9,9$ &\cr
height2pt&\omit&&\omit&&\omit&&\omit&\cr
& $5$  && $2^2\cdot 3\cdot 5\cdot 7^2$ && $5^2\cdot 117119$ && 
$1,1,10,10$  &\cr
height2pt&\omit&&\omit&&\omit&&\omit&\cr
} \hrule}}

\end{section}

\begin{section}{Harder's congruences}
Harder predicts a congruence between an elliptic modular form $f$ and a
Siegel modular form whenever a `large' prime $\ell$ divides a critical value
$L(f,s)$ of the $L$-series of the elliptic modular form, see \cite{Ha1,vdG1}. 
In cooperation with Harder we checked a few cases. This lends at the same
time credibility to our computations and conjectures and to Harder's conjecture.

For example, 
if $f$ is a newform in the $1$-dimensional space 
$S^{+}_{20}(\Gamma_0(2))$ then $61$ divides
$L(f,12)$ and one expects the congruence
$$
p^8+a(p)+p^{11} \equiv \lambda(p)  \, (\bmod \, 61)
$$
for the Hecke eigenvalues $\lambda(p)$ of an eigenform 
$F \in S_{2,10}(\Gamma_2[2])$ for every prime $p\neq 2$.
By using dimension formulas of R. Tsushima for $S_{j,k}(\Gamma_2(w^0))$ and 
$S_{j,k}(\Gamma_2(w^1))$ (see also Section~\ref{sec-Tsu}) 
we find that $\dim S_{2,10}(\Gamma_2(w^0))=0$ and 
$\dim S_{2,10}(\Gamma_2(w^1))=1$.
For $p \leq 37$ the eigenvalues $\lambda(p)$ we have calculated for a non-zero 
$F \in S_{2,10}(\Gamma_2(w^1))$
satisfy the required congruence, e.g.\ $\lambda(3)= 18360$
and $3^8-13092+3^{11}\equiv 18360  \, (\bmod \, 61)$. 

In the following table we list a few congruences that are 
valid for the eigenvalues that we find.
Also in these cases the corresponding spaces of modular forms 
are $1$-dimensional and the Siegel modular forms do not come from level~$1$.

\smallskip

\centerline{
\vbox{
\offinterlineskip
\hrule
\halign{&\vrule#& \quad \hfil#\hfil \strut \quad \cr
height2pt&\omit&&\omit&&\omit&&\omit& \cr
& $\langle f \rangle $  && $\langle F \rangle  $ && $s$  && $\ell$ & \cr
height2pt&\omit&&\omit&&\omit&&\omit&\cr
\noalign{\hrule}
height2pt&\omit&&\omit&&\omit&&\omit&\cr
& $S_{20}^{+}(\Gamma_0(2))$  && $S_{2,10}(\Gamma_2(w^1))$  && $12$ && $61$ &\cr
height2pt&\omit&&\omit&&\omit&&\omit&\cr
& $S_{20}^{+}(\Gamma_0(2))$  && $S_{10,6}(\Gamma_2(w^1))$  && $16$ && $109$ &\cr
height2pt&\omit&&\omit&&\omit&&\omit&\cr
& $S_{18}^{-}(\Gamma_0(2))$  && $S_{6,7}(\Gamma_2(w^1))$  && $13$ && $29$ &\cr
height2pt&\omit&&\omit&&\omit&&\omit&\cr
& $S_{20}^{-}(\Gamma_0(2))$  && $S_{12,5}(\Gamma_2(w^1))$  && $17$ && $79$ &\cr
height2pt&\omit&&\omit&&\omit&&\omit&\cr
& $S_{22}^{+}(\Gamma_0(2))$  && $S_{16,4}(\Gamma_2(w^1))$  && $20$ && $37$ &\cr
height2pt&\omit&&\omit&&\omit&&\omit&\cr
} \hrule}}

Moreover, for a newform $f \in S_{16}(\Gamma_0(4))$ we find that our traces of
Frobenius on both $S_{8,5}(\Gamma_2[2])^{[4,1^2]}$ and $S_{8,5}(\Gamma_2[2])^{[3^2]}$,
satisfy the expected congruence modulo $\ell=37$, which divides $L(f,13)$.

\end{section}

\end{document}